# APPROXIMATIONS OF THE WIENER SAUSAGE AND ITS CURVATURE MEASURES

By Jan Rataj[1], Evgeny Spodarev and Daniel Meschenmoser

*Charles University, Ulm University and Ulm University*

A parallel neighborhood of a path of a Brownian motion is sometimes called the *Wiener sausage*. We consider almost sure approximations of this random set by a sequence of random polyconvex sets and show that the convergence of the corresponding mean curvature measures holds under certain conditions in two and three dimensions. Based on these convergence results, the mean curvature measures of the Wiener sausage are calculated numerically by Monte Carlo simulations in two dimensions. The corresponding approximation formulae are given.

**1. Introduction.** Although the path of the Brownian motion is a random fractal set, the boundary of its parallel neighborhood of radius $r$ (also called the *Wiener sausage*) is a Lipschitz manifold and the closure of the complement to the parallel neighborhood has positive reach for almost all $r > 0$, at least in dimensions two and three; see [5]. Due to this fact, curvature measures of the Wiener sausage, including its volume, surface area and Euler–Poincaré characteristic, can be introduced, as in [13]. Since the 1930s, much attention has been paid to the computation of the mean volume of the Wiener sausage; see [8, 16] and [1]. Recently, its mean surface area was obtained in [12]. The same result was re-established, using different tools, in [9]. Explicit formulae for other mean curvature measures of the Wiener sausage are still unknown. The first attempt to treat them was made by Last (in [9]), who reduced the computation of his support measures to two *mean curvature functions*.

In this paper, we consider the curvature measures of the Wiener sausage in dimensions two and three. In the planar case, we show that the mean

Received November 2007; revised November 2008.

[1]Supported in part by Grants MSM 0021620839 and GACR 201/06/0302.

Euler number is finite, by using a general property of the set of critical values of the distance function and the self-similarity of the Brownian motion. Furthermore, we consider a.s. approximations of the Wiener sausage by polyconvex sets and show that the corresponding curvature measures converge almost surely and, in certain cases, also in mean. Finally, we attempt to study the mean curvatures of the Wiener sausages by simulation and find the corresponding approximation formulae. For that, we consider the approximating polyconvex sets whose curvature measures can be calculated numerically. Averaging over a large number of Monte Carlo realizations, the mean curvatures of these sets are estimated and the corresponding approximation formulae in two dimensions are found. Although the formulae for the mean area and boundary length are known explicitly (see the discussion in Section 7.1), they involve an integral of a functional of Bessel functions that has to be assessed numerically. Thus, apart from illustrating our convergence results stated in Sections 5 and 6, our approximation formulae for the mean area and boundary length are of practical interest. As for the mean Euler characteristic $\mathbb{E}V_0(\Xi_r)$ of the Wiener sausage, no analytical expression is known thus far. Thus, the approximation formula in Section 7.2 is novel and can be considered as a first step toward studying its properties.

The paper is organized as follows. In Section 2, curvature measures and intrinsic volumes are introduced for polyconvex sets and the sets of positive reach. Section 3 contains some preliminaries on Wiener sausages. In Section 4, we show that the mean Euler number of the Wiener sausage is finite in $\mathbb{R}^2$. The rest of the paper deals with approximations of the Wiener sausage. Namely, the convergence of parallel neighborhoods of compact sets in the Hausdorff metric is considered in Section 5. In Theorem 5.2, conditions are found under which the weak convergence of their curvature measures takes place. Polyconvex approximations of the Wiener sausage satisfying these conditions are described in Section 6. Finally, Section 7 contains the Monte Carlo simulation results for the mean curvature measures of the Wiener sausage.

**2. Intrinsic volumes and curvature measures.** Let $V_d$ be the Lebesgue measure and $\mathcal{H}^s$ the $s$-dimensional Hausdorff measure in $\mathbb{R}^d$ (see, e.g., [14]). Denote by $\mathcal{K}$ the set of all convex bodies in $\mathbb{R}^d$. Let $\mathcal{R}$ be the family of all finite unions of sets from $\mathcal{K}$. It is sometimes called the *convex ring* and its elements are referred to as *polyconvex sets*.

Let $A \oplus B$ be the pointwise sum of two sets, $A$ and $B$, in $\mathbb{R}^d$. For the closed ball $B = B_r(o)$ of radius $r \geq 0$ in $\mathbb{R}^d$ centered at the origin, the set $A_r = A \oplus B_r(o)$ is called the *r-parallel neighborhood of A*. The operation $A \mapsto A_r$ is known as *dilation*.

Let $A$ be a nonempty compact subset of $\mathbb{R}^d$. The *distance function* of $A$ is defined by

$$\Delta_A : x \mapsto \min\{|x - a| : a \in A\}, \qquad x \in \mathbb{R}^d.$$

Given two nonempty compact subsets $A, B$ of $\mathbb{R}^d$, their *Hausdorff distance* is defined as

$$d_H(A, B) = \max\Big\{\max_{a \in A} \Delta_B(a), \max_{b \in B} \Delta_A(b)\Big\}.$$

It is well known that $d_H$ is a metric. For any $x \in \mathbb{R}^d$, we denote by

$$\Sigma_A(x) = \{a \in A : |x - a| = \Delta_A(x)\}$$

the set of all points in $A$ which are nearest to $x$. The set $\Sigma_A(x)$ is always nonempty, by the compactness of $A$. We say that $x \in \mathbb{R}^d$ is a *critical point* of $\Delta_A$ if $x$ lies in the closed convex hull of $\Sigma_A(x)$. A point $x$ is called *regular* if it is not critical. A number $r > 0$ is a *critical value* of $\Delta_A$ is there exists a critical point $x$ of $\Delta_A$ with $\Delta_A(x) = r$. We shall denote by $C(A) \subseteq (0, \infty)$ the set of all critical values of $\Delta_A$.

The *reach* of a closed subset $A \subseteq \mathbb{R}^d$ is given by

$$\text{reach } A = \sup\{r \geq 0 : \forall x \in \mathbb{R}^d, \Delta_A(x) < r \Rightarrow \text{card } \Sigma_A(x) = 1\}$$

(see [3]). As examples of sets with positive reach, consider convex closed sets (with infinite reach) or sets with compact and $C^2$-smooth boundaries.

For a compact set $A \subseteq \mathbb{R}^d$ with positive reach, the Steiner formula holds for sufficiently small radii. More exactly, we have

$$V_d(A_r) = \sum_{i=0}^{d} \omega_i r^i V_{d-i}(A), \qquad 0 \leq r < \text{reach } A,$$

where $V_d$ is the Lebesgue measure in $\mathbb{R}^d$ (volume), $\omega_i$ is the volume of the unit $i$-ball and $V_i$ is the $i$th intrinsic volume of $A$, $i = 0, \ldots, d$ (see [3]). The functionals $V_i$ are motion invariant and additive. In particular, $V_0$ is the Euler–Poincaré characteristic and $V_{d-1}$ is one half of the surface area.

A local version of the Steiner formula holds for sets with positive reach. Let $\xi_A(x)$ be the nearest point of $A$ to $x$ whenever $\Delta_A(x) < \text{reach } A$. Then, we have, for any bounded Borel subset $F$ of $\mathbb{R}^d$,

$$V_d((A_r \setminus A) \cap \xi_A^{-1}(F)) = \sum_{i=1}^{d} \omega_i r^i C_{d-i}(A; F), \qquad 0 \leq r < \text{reach } A,$$

where $C_i(A; \cdot)$ is the $i$th *curvature measure* of $A$; it is a signed Radon measure concentrated on $\partial A$ for $0 \leq i \leq d - 1$ and $C_{d-1}(A; \cdot)$ is the restriction of the $(d - 1)$-dimensional Hausdorff measure to $\partial A$, provided that $A$ is $d$-dimensional. If $\partial A$ is compact, then $C_i(A; \mathbb{R}^d) = V_i(A)$, $i = 0, \ldots, d - 1$.

The curvature measures $C_i(\cdot;\cdot)$ admit an additive extension to the family of polyconvex sets; see [15]. The same holds, of course, for intrinsic volumes $V_i(\cdot)$, $i = 0, \ldots, d$. Curvature measures can further be extended to certain full-dimensional Lipschitz manifolds; see [13]. This applies, in particular, to the case when $A$ is a full-dimensional Lipschitz manifold in $\mathbb{R}^d$ and $\overline{\mathbb{R}^d \setminus A_r}$ has positive reach; see [13], Proposition 4. In this case, the curvature measures of $A_r$ can be introduced through the identity

$$C_i(A_r; F) = (-1)^{d-i-1} C_i(\overline{\mathbb{R}^d \setminus A_r}; F)$$

for $i = 0, \ldots, d-1$ and $r > 0$ such that $\text{reach}(\overline{\mathbb{R}^d \setminus A_r}) > 0$.

**3. Mean curvatures of the Wiener sausage.** Let $(\Omega, \mathfrak{F}, P)$ be an arbitrary probability space. The *Wiener process* is a random process $\{W(t) : t \geq 0\}$ defined on $(\Omega, \mathfrak{F}, P)$ with continuous paths and with the following properties:

- $W(0) = 0$ almost surely;
- $W$ has independent increments;
- $W(t) - W(s)$ has Gaussian distribution with mean zero and variance $t - s$ whenever $0 \leq s < t$.

Given $d$ independent Wiener processes $W_1, \ldots, W_d$, the random function $X(t) = (W_1(t), \ldots, W_d(t))$, $t \geq 0$, is the *Brownian motion* in $\mathbb{R}^d$ initiated at the origin $o \in \mathbb{R}^d$; see, for example, [2]. Denote by $\Xi = \{X(t) : 0 \leq t \leq 1\} \subseteq \mathbb{R}^d$ the image of the Brownian motion in $\mathbb{R}^d$ running over $0 \leq t \leq 1$.

DEFINITION 3.1. The set $\Xi_r = \Xi \oplus B_r(o)$, $r \geq 0$, is called a *Wiener sausage*; see, for example, [18], page 64.

By [12], Corollary 4.4, the intrinsic volumes $V_i(\Xi_r)$, $i = 0, \ldots, d-1$, are well defined for $d \leq 3$ almost surely. Moreover, the expectation $\mathbb{E}V_i(\Xi_r)$ is finite for $i = d, d-1$; see [1] and [12]. It will be shown in the next section that $\mathbb{E}V_0(\Xi_r) < \infty$ in the case $d = 2$. We conjecture that all mean intrinsic volumes of $\Xi_r$ are also finite in three dimensions.

**4. Finiteness of the mean Euler number.** In this section, we show that the mean Euler number $\mathbb{E}V_0(\Xi_r)$ exists and is finite for any $r > 0$ in two dimensions. Two basic tools are used. The first is a general result about the smallness of the set $C(A)$ of critical values of $A \subseteq \mathbb{R}^2$, which is of a similar nature to the main result in [5]. Furthermore, the stochastic self-similarity of a Brownian motion is used to show a certain local homogeneity of the random set $C(\Xi)$.

LEMMA 4.1. *Let* $A \subseteq \mathbb{R}^2$ *be compact with* $D = \text{diam}\, A$ *and let* $0 < \varepsilon < r$ *be such that* $C(A) \cap [r, r + \varepsilon) = \varnothing$. *Then, the total number of bounded connected components of* $\mathbb{R}^2 \setminus A_r$ *is less than or equal to* $\frac{\pi(D + 2r)^2}{r\sqrt{2r\varepsilon}}$.

PROOF. Let $L_1, L_2, \ldots$ be an enumeration of the bounded connected components of $\mathbb{R}^2 \setminus A_r$ and let $\varepsilon_i$ be the inradius of $L_i$ (i.e., $\varepsilon_i$ is the maximal number such that there exists a disc of radius $\varepsilon_i$ contained in $\overline{L_i}$). Let, for each $i$, $B_i = B_{\varepsilon_i}(s_i)$ be any fixed disc with this property. Then, $s_i$ is a critical point of $\triangle_A$ with critical value $r + \varepsilon_i$.

Consider two different centers, $s_i$ and $s_j$. Since they belong to different connected components, the segment $[s_i, s_j]$ must contain a point of $A_r$. This point must belong to a disc of radius $r$ which is disjoint with the interiors of both $B_i$ and $B_j$. Simple geometric reasoning leads to

$$|s_i - s_j| \geq \sqrt{2r\varepsilon_i + \varepsilon_i^2} + \sqrt{2r\varepsilon_j + \varepsilon_j^2}.$$

Let $Z_i = \{z \in \mathbb{R}^2 : \operatorname{dist}(z, B_i) \leq \operatorname{dist}(z, B_j) \forall j \neq i\}$ be the "zone of influence" of $B_i$. We shall find a lower bound for the area of $Z_i$. Let $b_i$ be a boundary point of $B_i$ belonging to $A_r$; thus, $a_i = b_i + \frac{r}{\varepsilon_i}(b_i - s_i) \in A$. Setting $s_i' = b_i + \frac{\varepsilon}{\varepsilon_i}(s_i - b_i)$, all three discs $B_i, B_r(a_i)$ and $B_\varepsilon(s_i')$ have a common tangent at $b_i$. Now, let $t_i^1$ and $t_i^2$ be two tangent points of $\partial B_r(a_i)$ on lines passing through $s_i'$. It follows from the previous observations that the quadrangle $a_i t_i^1 s_i' t_i^2$ lies in $Z_i$ and its area is $r\sqrt{2r\varepsilon + \varepsilon^2} \geq r\sqrt{2r\varepsilon}$. Consequently, the total number of zones $Z_i$ must not exceed $\frac{\pi(D+2r)^2}{r\sqrt{2r\varepsilon}}$. □

LEMMA 4.2. *Let $n \in \mathbb{N}$, $\varepsilon_1, \ldots, \varepsilon_n$ be positive numbers and $S_0, \ldots, S_n$ be segments of length at least $s > 0$ contained in the ball $B_R(o)$, $R \geq s$. Assume that for all $0 \leq i < j \leq n$, $\operatorname{dist}(S_i, S_j) \geq \sqrt{s/2}\sqrt{\varepsilon_{i+1} + \cdots + \varepsilon_j}$. Then, $\sqrt{\varepsilon_1} + \cdots + \sqrt{\varepsilon_n} \leq 8s^{-3/2}\pi R^2$.*

PROOF. To prove the lemma, we shall use the following fact which can be shown by induction. Whenever $n \in \mathbb{N}$, $\delta_1, \ldots, \delta_n > 0$ and $P_0, \ldots, P_n$ are points on a line such that

$$\operatorname{dist}(\{P_0, \ldots, P_{i-1}\}, \{P_i, \ldots, P_n\}) \geq \delta_i, \qquad i = 1, \ldots, n,$$

then $\delta_1 + \cdots + \delta_n \leq \operatorname{diam}\{P_0, \ldots, P_n\}$.

Now, let $L$ be any line hitting the disc $B_R(o)$ and let $I$ be the set of all indices $i \leq n$ such that $S_i$ hits $L$. Applying the fact to the points $P_i = L \cap S_i$ and $\delta_i = \sqrt{s/2}\sqrt{\varepsilon_i}$, $i \in I$, we get

$$\sqrt{\frac{s}{2}} \sum_{i \in I} \sqrt{\varepsilon_i} \leq \lambda^1(L \cap B_R(o)).$$

Let $\{e_1, e_2\}$ be an orthonormal basis of $\mathbb{R}^2$ and consider the two grids of lines $L_i^1 = e_1^\perp + (is/\sqrt{2})e_1$ and $L_i^2 = e_2^\perp + (is/\sqrt{2})e_2$, where $i$ are integers with $|i| \leq \sqrt{2}R/s$ and $e_j^\perp$ are the coordinate axes with normal vectors $e_j$, $j = 1, 2$. Note that each segment $S_i$ hits at least one of the two line grids. From the above fact, we get that

$$\sqrt{\frac{s}{2}} \sum_{i=1}^{n} \sqrt{\varepsilon_i} \leq 2 \sum_{|i| \leq \sqrt{2}R/s} \lambda^1(L_i^1 \cap B_R(o)) \leq 2(2R + \pi R^2/(s/\sqrt{2}))$$

$$\leq 2s^{-1}\pi R^2 \left( \frac{2}{\pi} \frac{s}{R} + \sqrt{2} \right) \leq 5s^{-1}\pi R^2,$$

hence

$$\sum_{i=1}^{n} \sqrt{\varepsilon_i} \leq 5\sqrt{2}s^{-3/2}\pi R^2 < 8s^{-3/2}\pi R^2. \qquad \square$$

LEMMA 4.3. *Let $A$ be a nonempty compact set in $\mathbb{R}^2$, $x, y$ two critical points of its distance function with distances $0 < r < s$ from $A$ and let $a, b$ be points of $A$ with $|a - x| = \Delta_A(x) = r$ and $|b - y| = \Delta_A(y) = s$. Then, $|x - y| \geq \sqrt{2r}\sqrt{s - r}$. Furthermore, if $S$ ($T$) is the segment with endpoints $x$ and $\frac{x+a}{2}$ (resp., $y$ and $\frac{y+b}{2}$), then $\mathrm{dist}(S, T) \geq \sqrt{r/2}\sqrt{s - r}$.*

PROOF. The interior of the ball $B_s(y)$ contains no point of $A$. Thus, $\mathrm{int}\, B_s(y)$ contains no diameter of $B_r(x)$ (since, otherwise, $x$ would not be a critical point). It follows that $|x - y| \geq \sqrt{2r}\sqrt{s - r}$ (see [5], Lemma 4.3). Furthermore, since $a \in \partial B_r(x) \setminus \mathrm{int}\, B_s(y)$ and $b \in \partial B_s(y) \setminus \mathrm{int}\, B_r(x)$, we have $(b - a) \cdot (y - x) \geq 0$. It follows that the orthogonal projections of the segments $S$ and $T$ onto the line passing through $x$ and $y$ have distance at least $|y - x|/2 \geq \sqrt{r/2}\sqrt{s - r}$. The same lower bound must then apply for the distance of the segments themselves. $\square$

COROLLARY 4.1. *Let $A \subseteq \mathbb{R}^2$ be compact with $D = \mathrm{diam}\, A$ and let $0 < s < r$. Write*

$$(s, r) \setminus C(A) = \bigcup_i I_i$$

*as a countable union of maximal disjoint open intervals with lengths $\varepsilon_i$. Then,*

$$\sum_i \sqrt{\varepsilon_i} \leq 8s^{-3/2}\pi D^2.$$

We now consider a Brownian motion. Given $r > 0$, we introduce the random variable

$$\eta_r = \inf\{\varepsilon > 0 : r + \varepsilon \in C(\Xi)\},$$

where we use the convention $\inf \varnothing = \infty$. By [12], Theorem 4.1, any $r > 0$ does not belong to $C(\Xi)$ almost surely. Hence, we get $\eta_r > 0$ almost surely. For any curve $\gamma : [0,1] \to \mathbb{R}^d$, define its supremum norm by $\|\gamma\| = \max_{t \in [0,1]} |\gamma(t)|$. When we write $\|X\|$ in the following, we implicitly restrict $X$ to $[0,1]$.

LEMMA 4.4. *Given $r > 0$, there exists a $\delta > 0$ such that for all $s \in (r - \delta, r)$ and all $\varepsilon < \delta$, we have*

$$P\left(\eta_s < \left(\frac{r}{r-\delta}\right)^3 \|X\|^4 \varepsilon\right) \geq \frac{1}{2} P(\eta_r < \|X\|^4 \varepsilon).$$

PROOF. For each trajectory of $X$, choose a critical point $z \in \mathbb{R}^2$ with $\Delta_\Xi(z) = r + \eta_r$ and denote the random variable $\xi = |X(1) - z| - (r + \eta_r)$. Since $X(1)$ does not lie on the boundary of $B_{r+\eta_r}(z)$ almost surely (see [12], Lemma 4.2 or [9], Lemma 3.3), we have $\xi > 0$ a.s. Thus, as $\tau \searrow 0$, it holds that

$$\{\eta_r < \|X\|^4 \varepsilon, \xi > \tau\} \nearrow \{\eta_r < \|X\|^4 \varepsilon\}.$$

By choosing $\tau > 0$ small enough, we have

$$P(\eta_r < \|X\|^4 \varepsilon, \xi > \tau) \geq \tfrac{3}{4} P(\eta_r < \|X\|^4 \varepsilon).$$

Furthermore, we can choose $\lambda_0 > 1$ such that

$$P\left(\sup_{1 < t < \lambda_0} |X(t) - X(1)| < \tau\right) \geq \frac{3}{4}.$$

Since a Brownian motion has independent increments, we have

$$P\left(\eta_r < \|X\|^4 \varepsilon, \xi > \tau, \sup_{1 < t < \lambda_0} |X(t) - X(1)| < \tau\right)$$

$$= P(\eta_r < \|X\|^4 \varepsilon, \xi > \tau) P\left(\sup_{1 < t < \lambda_0} |X(t) - X(1)| < \tau\right)$$

$$\geq \frac{9}{16} P(\eta_r < \|X\|^4 \varepsilon).$$

Note that if $\eta_r < \varepsilon$, $\xi > \tau$ and $\sup_{1 < t < \lambda_0} |X(t) - X(1)| < \tau$, then $r + \eta_r$ remains a critical value of $\{X(t) : 0 \leq t \leq \lambda_0\}$. By the scaling invariance property, $\Xi^* = \lambda^{-1/2}\{X(t) : 0 \leq t \leq \lambda\}$ has the same distribution as $\Xi$, with critical value $\lambda^{-1/2}(r + \eta_r)$ for any $1 \leq \lambda \leq \lambda_0$. Denoting by $\|X\|^*, \eta_s^*$ the random variables corresponding to $\Xi^*$ [thus, $\|X\|^* = \lambda^{-1/2}\max_{0 \leq t \leq \lambda} |X(t)|$], we have $\eta_{\lambda^{-1/2}r}^* \leq \lambda^{-1/2}\eta_r$, $\|X\|^* \geq \lambda^{-1/2}\|X\|$. Consequently, we have

$$\frac{1}{2} P(\eta_r < \|X\|^4 \varepsilon) < \frac{9}{16} P(\eta_r < \|X\|^4 \varepsilon)$$

$$\leq P\left(\eta_r < \|X\|^4 \varepsilon, \xi > \tau, \sup_{1 < t < \lambda_0} |X(t) - X(1)| < \tau\right)$$

$$= P\left(\eta^*_{\lambda^{-1/2}r} < \lambda^{-1/2}\|X\|^4\varepsilon, \xi > \tau, \sup_{1 < t < \lambda_0} |X(t) - X(1)| < \tau\right)$$

$$\leq P\left(\eta^*_{\lambda^{-1/2}r} < \lambda^{3/2}(\|X\|^*)^4\varepsilon, \xi > \tau, \sup_{1 < t < \lambda_0} |X(t) - X(1)| < \tau\right)$$

$$= P\left(\eta_{\lambda^{-1/2}r} < \lambda^{3/2}\|X\|^4\varepsilon, \xi > \tau, \sup_{1 < t < \lambda_0} |X(t) - X(1)| < \tau\right)$$

$$\leq P(\eta_{\lambda^{-1/2}r} < \lambda^{3/2}\|X\|^4\varepsilon)$$

for all $1 \leq \lambda \leq \lambda_0$. The assertion thus holds with $\delta = r - \lambda_0^{-1/2}r$.   $\square$

PROPOSITION 4.1.   *For any $r > 0$, the mean number of connected components of $\mathbb{R}^2 \setminus \Xi_r$ is finite.*

PROOF.   Let $N_r$ be the number of bounded connected components of $\mathbb{R}^2 \setminus \Xi_r$. For $r > 0$, consider $\eta_r$. If $\eta_r = \infty$, then $\Xi_r$ does not have any holes and hence $N_r = 0$ [otherwise, any hole would contain a maximal inscribed circle with its center at a critical point and with critical radius $s \in (r, \infty)$]. Since $r \notin C(\Xi)$ a.s., we get $C(\Xi) \cap [r, r + \eta_r) = \varnothing$ a.s. By Lemma 4.1, we thus have

$$N_r \leq \frac{\pi}{r\sqrt{2r}} \frac{(2\|X\| + 2r)^2}{\sqrt{\eta_r}} \qquad \text{a.s.},$$

where we adopt the convention that $\frac{1}{\infty} = 0$. We shall show that $\mathbb{E}\|X\|^2\eta_r^{-1/2}$ is finite, which is equivalent to the assertion. Let $\sigma = \frac{r}{r-\delta}$. We have

$$\mathbb{E}\|X\|^2\eta_r^{-1/2} = \int_0^\infty P(\|X\|^2\eta_r^{-1/2} > t)\,dt = \int_0^\infty P(\eta_r < \|X\|^4 t^{-2})\,dt$$

$$\leq \frac{2}{\delta}\int_0^\infty \mathbb{E}\int_{r-\delta}^r \mathbf{1}_{(\eta_s < \sigma^3\|X\|^4 t^{-2})}\,ds\,dt,$$

by Lemma 4.4 and Fubini's theorem, since

$$\frac{1}{2}\delta P(\eta_r < \|X\|^4\varepsilon) = \frac{1}{2}\int_{r-\delta}^r P(\eta_r < \|X\|^4\varepsilon)\,ds$$

$$\leq \int_{r-\delta}^r P(\eta_s < \sigma^3\|X\|^4\varepsilon)\,ds$$

$$= \mathbb{E}\int_{r-\delta}^r \mathbf{1}_{(\eta_s < \sigma^3\|X\|^4\varepsilon)}\,ds$$

for any $\varepsilon > 0$. Fix a trajectory $\Xi$ of $X$. Let $\{I_i\}$ be the family of maximal open intervals of the complement of $C(\Xi)$ which intersect the interval $(r - \delta, r)$ and let $\varepsilon_i$ be the length of $I_i$. We assume, without loss of generality, that the sequence $\{\varepsilon_i\}$ is nonincreasing. We then have

$$\int_{r-\delta}^r \mathbf{1}_{(\eta_s < \sigma^3 \|X\|^4 t^{-2})} \, ds \leq \sum_i \min\{\sigma^3 \|X\|^4 t^{-2}, \varepsilon_i\}$$

and hence, by Corollary 4.1,

$$\mathbb{E}\|X\|^2 \eta_r^{-1/2} \leq \frac{2}{\delta} \mathbb{E} \int_0^\infty \sum_i \min\{\sigma^3 \|X\|^4 t^{-2}, \varepsilon_i\} \, dt$$

$$= \frac{2}{\delta} \mathbb{E} \sum_i \left( \int_0^{\sigma^{3/2} \|X\|^2 \varepsilon_i^{-1/2}} \varepsilon_i \, dt + \int_{\sigma^{3/2} \|X\|^2 \varepsilon_i^{-1/2}}^\infty t^{-2} \, dt \right)$$

$$= \frac{2}{\delta} \mathbb{E} \sum_i (\sigma^{3/2} \|X\|^2 \sqrt{\varepsilon_i} + \sigma^{-3/2} \|X\|^{-2} \sqrt{\varepsilon_i})$$

$$\leq \frac{2}{\delta} \mathbb{E}(\sigma^{3/2} \|X\|^2 + \sigma^{-3/2} \|X\|^{-2}) 8(r-\delta)^{-3/2} \pi \|X\|^2$$

$$= \frac{16\pi}{\delta(r-\delta)^{3/2}} (\sigma^{3/2} \mathbb{E}\|X\|^4 + \sigma^{-3/2}).$$

Since all moments of $\|\Xi\|$ are finite, the proposition is proved. $\square$

COROLLARY 4.2. *The mean Euler number of the Wiener sausage $\Xi_r$ in $\mathbb{R}^2$ is finite for any $r > 0$.*

**5. Convergence of parallel neighborhoods.** Let $A \subseteq \mathbb{R}^d$ be a nonempty and compact subset of $\mathbb{R}^d$. For $r \geq 0$, denote by $V_A(r)$ the volume $V_d(A_r)$ of the parallel neighborhood to $A$. The following result was proven by Stachó [17], Theorem 3.

THEOREM 5.1. *Let $\{A^n\}$ be a sequence of nonempty compact subsets of $\mathbb{R}^d$ converging with respect to the Hausdorff metric to a nonempty compact set $A$. Then,*

(i) *$V_{A^n}(r) \to V_A(r)$, $n \to \infty$, for all $r > 0$;*
(ii) *$\mathcal{H}^{d-1}(\partial A_r^n) \to \mathcal{H}^{d-1}(\partial A_r)$, $n \to \infty$, for all $r > 0$ where $(V_A)'(r)$ exists.*

REMARK 5.1. In fact, assertion (i) is included in the proof of [17], Theorem 3. It follows from the fact that for any $\varepsilon > 0$, $A_{r-\varepsilon} \subseteq (A^n)_r \subseteq A_{r+\varepsilon}$ for sufficiently large $n$ and from the continuity of the volume function $V_A$. Assertion (ii) is formulated in a different form in [17], Theorem 3, but (ii) is exactly what is proved there.

We shall now show that the curvature measures (and, hence, also the intrinsic volumes) of the parallel neighborhood also converge.

THEOREM 5.2.    *Let $\{A^n\}$ be a sequence of nonempty compact subsets of $\mathbb{R}^d$ converging with respect to the Hausdorff metric to a nonempty compact set $A$ and let $r \in (0, \infty) \setminus C(A)$. Then, $r \notin C(A^n)$ for sufficiently large $n$ and the curvature measures $C_i(A_r^n; \cdot)$ converge weakly to $C_i(A_r; \cdot)$, $n \to \infty$, for $i = 0, \ldots, d-1$.*

The proof will be based on the following result.

THEOREM 5.3 ([3], Theorem 5.9).    *Let $D^n, D$ be compact sets with positive reach such that $d_H(D^n, D) \to 0$ and $\inf_n \operatorname{reach} D^n > 0$. Then, $C_i(D^n; \cdot) \to C_i(D; \cdot)$ weakly, $n \to \infty$, $i = 0, \ldots, d-1$.*

To prove Theorem 5.2, it is sufficient to show that

$$\liminf_{n \to \infty} \operatorname{reach} C_r^n > 0,$$

where $C_r^n = \overline{\mathbb{R}^d \setminus A_r^n}$. We shall start with some continuity results.

We say that a function $F$ from a metric space $X$ to the space $\mathcal{K}'$ of nonempty compact subsets of $\mathbb{R}^d$ (equipped with the Hausdorff metric) is *upper semicontinuous* if, whenever $x_n \to x$ in $X$, $a_n \to a$ in $\mathbb{R}^d$ and $a_n \in F(x_n)$, $n \in \mathbb{N}$, we have $a \in F(x)$ (see [11]). Recall, also, that a real function $f : X \to \mathbb{R}$ is upper semicontinuous if $\liminf_{y \to x} f(y) \leq f(x)$ for all $x \in X$.

LEMMA 5.1.    (i) *The function $x \mapsto \Sigma_A(x)$ is upper semicontinuous on $\mathbb{R}^d$ for any compact $A \subseteq \mathbb{R}^d$.*

(ii) *The function $A \mapsto \Sigma_A(x)$ is upper semicontinuous on $\mathcal{K}'$ for any $x \in \mathbb{R}^d$.*

PROOF.    To verify (i), let $x_n \to x$, $a_n \in \Sigma_A(x_n)$, $a_n \to a$. We shall show that $a \in \Sigma_A(x)$. If not, there would be another point $b \in A$ with $|b - x| < |a - x|$. Let $n$ be sufficiently large that

$$\max\{|a_n - a|, |x_n - x|\} < \varepsilon = \frac{1}{3}(|a - x| - |b - x|).$$

Then, by the triangle inequality,

$$|x_n - b| \leq |x_n - x| + |b - x| < \varepsilon + |a - x| - 3\varepsilon$$
$$< |a - x| - |a - a_n| - \varepsilon \leq |a_n - x| - \varepsilon$$
$$< |a_n - x| - |x - x_n| \leq |a_n - x_n|,$$

which means that $a_n$ is not the closest point of $A$ to $x_n$, a contradiction.

In order to show (ii), let $A^n, A$ be compact sets, $d_H(A^n, A) \to 0$, $a_n \in \Sigma_{A^n}(x)$, $a_n \to a$. We have to show that $a \in \Sigma_A(x)$. If not, there would be a point $b \in A$ with $|b - x| < |a - x|$. By the definition of Hausdorff distance, there is a point $b_n \in A_n$ with $|b_n - b| < \varepsilon$ for sufficiently large $n$, where $\varepsilon = \frac{1}{2}(|a - x| - |b - x|)$. We also have $|a_n - a| < \varepsilon$ for large $n$. We then easily obtain, by the triangle inequality, that $|x - b_n| < |x - a_n|$, which contradicts the fact that $a_n$ is the closest point of $A_n$ to $x$. $\quad\square$

As in [12], we introduce the function

$$J_A : x \mapsto \min\{|a - x| : a \in \operatorname{conv} \Sigma_A(x)\}, \qquad x \in \mathbb{R}^d.$$

Clearly, $x$ is a regular point of $\Delta_A$ if and only if $J_A(x) > 0$.

LEMMA 5.2.  (i) *The function* $x \mapsto J_A(x)$ *is upper semicontinuous on* $\mathbb{R}^d$ *for any compact* $A \subseteq \mathbb{R}^d$.

(ii)  *The function* $A \mapsto J_A(x)$ *is upper semicontinuous on* $\mathcal{K}'$ *for any* $x \in \mathbb{R}^d$.

PROOF.  This follows directly from Lemma 5.1 and from the fact that $J_A(x)$ depends continuously on $\Sigma_A(x)$. $\quad\square$

LEMMA 5.3.  *If* $A^n, A$ *are nonempty compact subsets of* $\mathbb{R}^d$, $d_H(A^n, A) \to 0$ *as* $n \to \infty$ *and* $r \notin C(A)$, *then*

$$\liminf_{n \to \infty} \inf\{J_{A^n}(x) : x \in \partial A_r^n\} > 0.$$

PROOF.  We have $\inf\{J_A(x) : x \in \partial A_r\} > 0$ by [12], Lemma 3.3. Further, since $J_A$ is upper semicontinuous and $\partial A$ is compact, we have

$$\eta = \inf\{J_A(x) : r - \varepsilon \le \Delta_A(x) \le r + \varepsilon\} > 0$$

for some sufficiently small $\varepsilon > 0$. By the definition of the Hausdorff metric,

$$A_{r-\varepsilon} \subseteq A_r^n \subseteq A_{r+\varepsilon}$$

for large $n$. Now, for any $x$ with $r - \varepsilon \le \Delta_A(x) \le r + \varepsilon$, we have $J_{A^n}(x) > \eta/2$ for sufficiently large $n$, by Lemma 5.2. Define

$$K^n = \{x : r - \varepsilon \le \Delta_A(x) \le r + \varepsilon, J_{A^m}(x) \le \eta/2 \text{ for } m \ge n\}.$$

By the upper semicontinuity of $J_{A^m}$, the sets $K_n$ are compact and since $K^n \searrow \varnothing$, we have $K^{n_0} = \varnothing$ for some $n_0$, which completes the proof. $\quad\square$

LEMMA 5.4.  *If* $A \subseteq \mathbb{R}^d$ *is a nonempty and compact subset of* $\mathbb{R}^d$ *and* $r > 0$, *then*

$$\operatorname{reach} \overline{\mathbb{R}^d \setminus A_r} \ge \inf\{J_A(x) : x \in \partial A_r\}.$$

PROOF.    If $z \in \partial A_r$, then, for any $a \in \operatorname{conv} \Sigma_A(z)$, $B_{|z-a|}(a) \subseteq A_r$ with

$$|z - a| \geq \eta = \inf\{J_A(x) : x \in \partial A_r\}.$$

We can thus "roll" a ball of radius $\eta$ from outside along the boundary of $\overline{\mathbb{R}^d \setminus A_r}$, which proves that reach $\overline{\mathbb{R}^d \setminus A_r} \geq \eta$ (see [3], Theorem 4.18).    □

PROOF OF THEOREM 5.2.    Follows directly from Lemmas 5.3 and 5.4.    □

**6. Approximations of the Wiener sausage by polyconvex sets.**  In this section, we use a piecewise linear approximation $\Xi^n$ of the path $\Xi$ of a Brownian motion $X$ to construct the almost sure approximation of the Wiener sausage $\Xi_r = \Xi \oplus B_r(o)$ by random polyconvex sets $\Xi_r^n = \Xi^n \oplus B_r(o)$. Under some assumptions, the convergence of the corresponding curvature measures and their first moments is shown in Corollary 6.1 and Proposition 6.1. In Sections 6.2 and 6.3, examples of such approximations are considered in detail.

6.1. *General convergence results.*  Let $\Xi^n$ be an arbitrary random piecewise linear curve with vertices $X^n(t_i)$, $i = 1, \ldots, k_n$, approximating the path $\Xi$ of the Brownian motion in the sense of the Hausdorff distance: $d_H(\Xi^n, \Xi) \to 0$ as $n \to \infty$, almost surely. For each $n \in \mathbb{N}$, consider its parallel neighborhood $\Xi_r^n = \Xi^n \oplus B_r(o)$, $r > 0$. It holds that $\Xi_r^n \in \mathcal{R}$ a.s. for all $n \in \mathbb{N}$ and $r > 0$. We call the sequence $\{\Xi_r^n\}_{n \in \mathbb{N}}$ of polyconvex random closed sets an *(almost sure) polyconvex approximation* of the Wiener sausage $\Xi_r$.

Since $\Xi_r^n$ is polyconvex, its curvature measures $C_i(\Xi_r^n; \cdot)$ and intrinsic volumes $V_i(\Xi_r^n)$, $i = 0, \ldots, d$, are well defined for all $r > 0$. In the following, we show that the curvature measures of $\Xi_r^n$ approximate the corresponding curvature measures of the Wiener sausage for growing $n$.

The next result follows immediately from Theorems 5.1, 5.2 and [12], Corollary 4.1 and Lemma 4.3.

COROLLARY 6.1.    *Let $d_H(\Xi^n, \Xi) \to 0$ as $n \to \infty$, almost surely. We then have the following:*

(i) *$V_{\Xi^n}(r) \to V_\Xi(r)$ as $n \to \infty$ almost surely for any $r > 0$;*
(ii) *$\mathcal{H}^{d-1}(\partial \Xi_r^n) \to \mathcal{H}^{d-1}(\partial \Xi_r)$ as $n \to \infty$ almost surely for any $r > 0$ if $d \leq 3$, and for almost all $r > 0$ if $d \geq 4$;*
(iii) *$C_i(\Xi_r^n; \cdot) \to C_i(\Xi_r; )$ weakly as $n \to \infty$ almost surely for any $r > 0$ and $i = 0, \ldots, d-1$ if $d \leq 3$.*

Corollary 6.1 combined with the dominated convergence theorem yields the following proposition.

PROPOSITION 6.1. *Let $d_H(\Xi^n, \Xi) \to 0$ as $n \to \infty$, almost surely. We then have the following:*

(i) $\mathbb{E} V_{\Xi^n}(r) \to \mathbb{E} V_{\Xi}(r)$ *as* $n \to \infty$ *for any* $r > 0$;

(ii) $\mathbb{E} \mathcal{H}^{d-1}(\partial \Xi^n_r) \to \mathbb{E} \mathcal{H}^{d-1}(\partial \Xi_r)$ *as* $n \to \infty$ *for any* $r > 0$ *if* $d \le 3$, *and for almost all* $r > 0$ *if* $d \ge 4$.

PROOF. It suffices to find the uniform integrable random upper bounds for $|V_i(\Xi^n_r)|$, $i = d, d-1$.

(i) First, note that

$$|\|X^n\| - \|X\|| \le d_H(\Xi^n, \Xi).$$

Since the last expression tends to 0 as $n \to 0$, we have

$$V_d(\Xi^n_r) \le V_d(B_{r+\|X^n\|}(o)) = \omega_d(r + \|X^n\|)^d \le \omega_d(r + \|X\| + 1)^d$$

for sufficiently large $n$ (and all $\omega$). It is known that all moments of the maximum of the Bessel process $|X(t)|$ are finite; see, for example, [12], proof of Theorem 2.1. Consequently, we have

$$\mathbb{E} V_d(\Xi^n_r) \le \omega_d \mathbb{E}(r + \|X\| + 1)^d < \infty$$

for all sufficiently large $n$.

(ii) Since $\partial \Xi^n_r \subseteq \Delta^{-1}_{\Xi^n}(\{r\})$ holds for all $r > 0$, we use the co-area formula from [4] to obtain

$$\mathcal{H}^{d-1}(\Xi^n_r) \le \mathcal{H}^{d-1}(\Delta^{-1}_{\Xi^n}(\{r\})) = V'_{\Xi_n}(r)$$

for all $r > 0$ such that $V'_{\Xi_n}(r)$ exists. These are all $r \notin C(\Xi_n)$. By Theorem 5.2, it holds that $r \notin C(\Xi_n)$ for all sufficiently large $n$ if $r \notin C(\Xi)$. It is known (see [12]) that $r \notin C(\Xi)$ almost surely for any $r > 0$ if $d \le 3$, and for almost all $r > 0$ if $d \ge 4$. By [12], Lemma 4.6, we obtain

$$\mathcal{H}^{d-1}(\Xi^n_r) \le V'_{\Xi_n}(r) \le \frac{d}{r} \omega_d(R + r)^d,$$

where $R = \|X^n\| \le \|X\| + \varepsilon$ for all sufficiently large $n$. Since the upper bound

$$\frac{d}{r} \omega_d(\|X\| + r + \varepsilon)^d$$

is integrable, the dominated convergence theorem yields $\mathbb{E} \mathcal{H}^{d-1}(\partial \Xi^n_r) \to \mathbb{E} \mathcal{H}^{d-1}(\partial \Xi_r)$ as $n \to \infty$. $\square$

The following statement follows directly from Corollary 6.1(iii) and the dominated convergence theorem.

REMARK 6.1.   Let $d_H(\Xi^n, \Xi) \to 0$ as $n \to \infty$, almost surely, $d \le 3$. If the sequence of random variables $\{V_i(\Xi_r^n)\}_{n \in \mathbb{N}}$ is uniformly integrable for $i = 0, \ldots, d-2$, then $\mathbb{E}V_i(\Xi_r^n) \to \mathbb{E}V_i(\Xi_r)$ as $n \to \infty$ for $i = 0, \ldots, d-2$ and $r > 0$. We conjecture that in the two-dimensional case, the sequence $\{V_0(\Xi_r^n)\}_{n \in \mathbb{N}}$ is uniformly integrable and that convergence for the mean Euler number holds.

6.2. *An example.*   Consider a regular lattice of time moments $t_i = i/k_n$, $i = 1, \ldots, k_n$, $k_n \to \infty$ as $n \to \infty$. Let $X_n$ be the piecewise linear approximation of $X$ passing through $X^n(t_i) = X(t_i)$ for all $i = 1, \ldots, k_n$ and $n \in \mathbb{N}$. It means that the vertices of the piecewise linear approximation $X^n$ of the Brownian motion $X$ lie on the path $\Xi^n$ of the Brownian motion. In the following, we show that such an approximation procedure preserves the convergence of the mean curvatures in the sense of the last section.

THEOREM 6.1.   *It holds that*

$$d_H(\Xi^n, \Xi) \le \|X^n - X\|, \qquad n \in \mathbb{N},$$

*almost surely.*

PROOF.   Using the piecewise linearity of $\Xi^n$, we get

$$\max_{y \in \Xi} \Delta_{\Xi_n}(y) = \max_{i=1,\ldots,k_n} \max_{t \in [(i-1)/k_n, i/k_n]} \Delta_{\Xi_n}(X(t))$$

$$= \max_{i=1,\ldots,k_n} \max_{t \in [(i-1)/k_n, i/k_n]} \min_{i=1,\ldots,k_n} \min_{s \in [(i-1)/k_n, i/k_n]} |X^n(s) - X(t)|$$

$$\le \max_{i=1,\ldots,k_n} \max_{t \in [(i-1)/k_n, i/k_n]} |X^n(t) - X(t)| = \max_{t \in [0,1]} |X^n(t) - X(t)|.$$

Using similar arguments, it can be shown that

$$\max_{y \in \Xi_n} \Delta_{\Xi}(y) \le \|X^n - X\|.$$

By the definition of Hausdorff distance, we have

$$d_H(\Xi^n, \Xi) = \max\left\{ \max_{y \in \Xi^n} \Delta_{\Xi}(y), \max_{y \in \Xi} \Delta_{\Xi^n}(y) \right\} \le \|X^n - X\|. \qquad \square$$

COROLLARY 6.2.   $d_H(\Xi^n, \Xi) \to 0$ *as* $n \to \infty$, *almost surely.*

PROOF.   It is known that $\|W_i^n - W_i\| \to 0$ as $n \to \infty$ almost surely for all $i = 1, \ldots, d$. The assertion of the corollary follows from the above relation, Theorem 6.1 and the inequality

$$\|X^n - X\| = \max_{t \in [0,1]} \sqrt{\sum_{i=1}^d (W_i^n(t) - W_i(t))^2} \le \sum_{i=1}^d \|W_i^n - W_i\| \to 0. \qquad \square$$

6.3. *Haar–Schauder approximation.* The piecewise linear approximations considered above can easily be simulated using Gaussian finite-dimensional distributions of a Brownian motion. Another constructive approach is provided by the so-called *Haar–Schauder approximation* of a Brownian motion; see, for example, the book [6].

Let $H_k : [0,1] \to \mathbb{R}$, $k \in \mathbb{N}$, be the so-called *Haar* function defined by the relations

$$H_1(s) = 1, \qquad s \in [0,1],$$

$$H_{2^m+k}(s) = \begin{cases} 2^{m/2}, & s \in \left[ \dfrac{k-1}{2^m}, \dfrac{2k-1}{2^{m+1}} \right), \\[2mm] -2^{m/2}, & s \in \left[ \dfrac{2k-1}{2^{m+1}}, \dfrac{k}{2^m} \right), \\[2mm] 0, & \text{otherwise}, \end{cases}$$

for $k = 1, 2, \ldots, 2^m$ and $m = 0, 1, 2, \ldots$. The *Schauder function* is given by

$$S_k(t) = \int_0^t H_k(s) \, ds, \qquad k \in \mathbb{N}.$$

It is known that the Wiener process $W$ is equal in distribution to the *Haar–Schauder series*

$$W(t) = \sum_{k=1}^\infty Y_k S_k(t), \qquad t \in [0,1],$$

where $\{Y_n\}_{n \in \mathbb{N}}$ is the sequence of i.i.d. $N(0,1)$-distributed random variables; see [6], pages 56–59. This series converges (a.s.) absolutely and uniformly on $[0,1]$. Hence, the Wiener process can be approximated pathwise by partial sums

$$\sum_{k=1}^{2^n} Y_k S_k(t), \qquad t \in [0,1].$$

We use this idea to approximate the coordinates $W_i(t)$, $i = 1, \ldots, d$, of the Brownian motion $X(t)$ a.s. by

$$W_i^n(t) = \sum_{k=1}^{2^n} Y_{ik} S_k(t), \qquad t \in [0,1],$$

where the sequences $\{Y_{ik}\}_{k \in \mathbb{N}}$, $i = 1, \ldots, d$, of i.i.d. standard normally distributed random variables $Y_{ik}$ are independent. For

$$X^n(t) = (W_1^n(t), \ldots, W_d^n(t)), \qquad t \in [0,1],$$

consider its path $\Xi^n = \{X^n(t) : 0 \le t \le 1\}$. It is a piecewise linear curve with $2^n + 1$ nodes lying on $\Xi$. More precisely, we have $X_n(0) = X(0) = o$ a.s. for

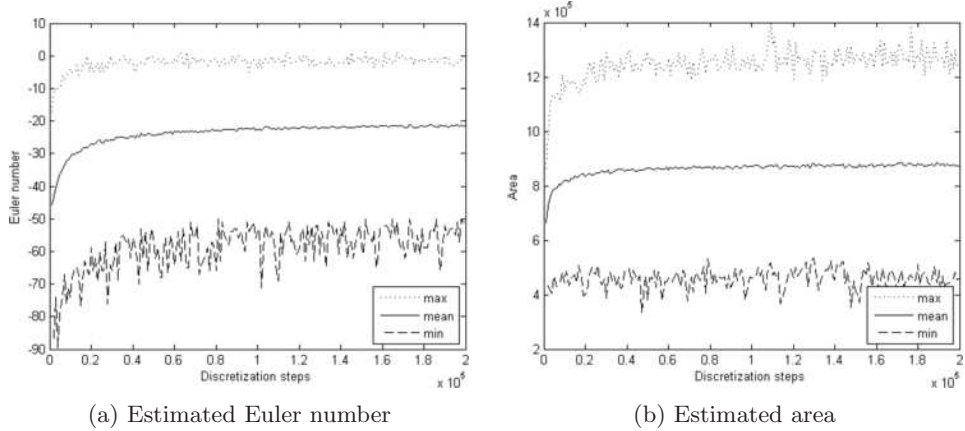

(a) Estimated Euler number                (b) Estimated area

Fig. 1.   *Mean (solid), minimal (dashed) and maximal (dotted) estimated Euler number and area, respectively, of 1000 simulated Wiener sausages, depending on the number of discretization steps.*

$n \in \mathbb{N}$ and

$$X_n(t) = (1 - \alpha_n(t))X\left(\frac{k-1}{2^n}\right) + \alpha_n(t)X\left(\frac{k}{2^n}\right), \qquad t \in \left[\frac{k-1}{2^n}, \frac{k}{2^n}\right]$$

a.s. for all $n$ and $k = 1, \ldots, 2^n$, where $\alpha_n(t) = 2^n t - k + 1$. Thus, $\Xi^n$ can be regarded as an a.s. approximation of the path $\Xi$ by piecewise linear curves, in the sense of Section 6.3, with $k_n = 2^n$.

**7. Numerical results.**   In this section, we present some numerical results for the estimated mean intrinsic volumes of the Wiener sausage in the two-dimensional case. The theoretical basis for the simulation which supports the outcome is provided by the convergence results of Proposition 6.1 and Corollary 6.2.

7.1. *Approach and estimation results.*   A piecewise linear approximation of the path of a Brownian motion in $\mathbb{R}^2$ is considered and the path is dilated to produce a polyconvex approximation of the Wiener sausage. The approximation of the path is achieved by simulating the independent, normally distributed increments of the Wiener process; see Section 6.2. For this, we chose the time interval $[0, T] = [0, 10^6]$ and $k$ equidistant time points in this interval. Figure 1 shows the estimated Euler number and estimated area, respectively, of approximated Wiener sausages $\Xi_r$ in the plane with radius $r = 20$ depending on the number $k$ of such discretization points. For each $k$, one can see the mean (solid) as well as the minimal (dashed) and maximal (dotted) values of 1000 realizations. It is obvious that for small values $k$, the approximation is not very good and the estimated value of the polyconvex approximation is therefore too small. On the other hand, both estimated values converge to a limit once $k$ is sufficiently large. Since this convergence behavior is basically the same for the boundary length, we discretized the interval $[0, T]$ in $k = 100{,}001$ time points for all simulations. The radius $r$ of the Wiener sausage varies between 1 and 1900 pixels. For radii less than 1, there is no difference between the Brownian motion and the Wiener sausage in the simulation due to discretization. In any case, in Figures 2(a)–2(f), one can recognize a trend in the curves for small radii. On the other hand, the behavior of the curves for growing radii is clearly visible, so there is no need to simulate bigger radii. In Figures 2(a)–2(f), the solid line denotes the mean value, the dashed line the minimal value and the dotted line the maximal value of the considered estimator out of 1000 realizations for each radius. The intrinsic volumes were estimated with the algorithm described in [7] using Steiner's formula. Figure 2(a) shows the mean, minimal and maximal estimated area of 1000 realizations. For growing radii, the area grows quadratically and as the radius tends to zero, the area tends to zero as well. This can be seen in Figure 2(b), which shows a magnification of Figure 2(a) for small radii. Thus, we get an empirical illustration of the following asymptotic result of Le Gall [10]:

$$\mathbb{E}V_{\Xi}(r) \sim \frac{\pi T}{|\log r|}, \qquad r \to 0.$$

Figures 2(c) and 2(d) show the estimated mean boundary length. One can see that it grows linearly once all holes of $\Xi_r$ disappear as $r$ tends to

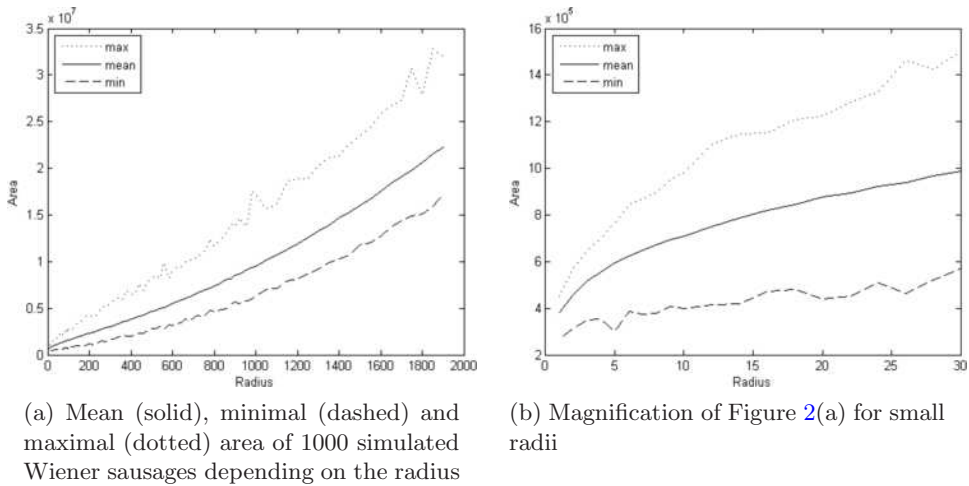

(a) Mean (solid), minimal (dashed) and maximal (dotted) area of 1000 simulated Wiener sausages depending on the radius

(b) Magnification of Figure 2(a) for small radii

FIG. 2. *Estimated mean intrinsic volumes of the Wiener sausage in two dimensions depending on the radius.*

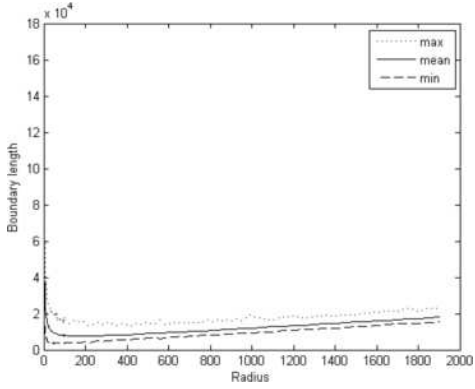

(c) Mean (solid), minimal (dashed) and maximal (dotted) boundary length of 1000 simulated Wiener sausages depending on the radius

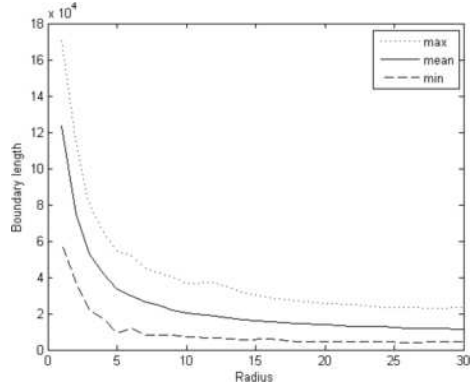

(d) Magnification of Figure 2(c) for small radii

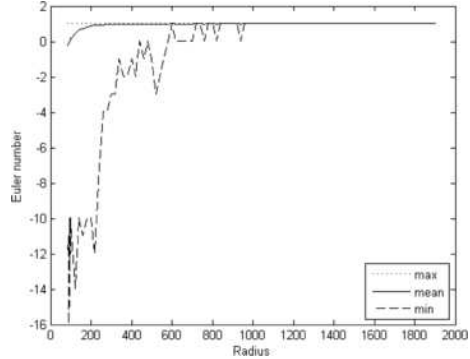

(e) Mean (solid), minimal (dashed) and maximal (dotted) Euler number of 1000 simulated Wiener sausages with large radii ($r \geq 120$)

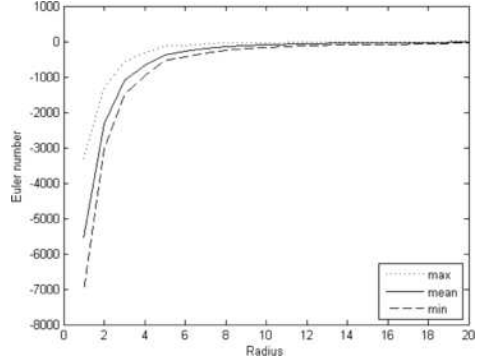

(f) Mean (solid), minimal (dashed) and maximal (dotted) Euler number of 1000 simulated Wiener sausages with small radii ($r \leq 15$)

Fig. 2. *Continued.*

infinity. It is shown in [12] that the expected boundary length equals the first derivative of the mean area of the Wiener sausage with respect to its radius. This is also backed by our simulations. For instance, the gradient of the tangent line to the graph of the mean area in Figure 2(b) grows rapidly as the radius becomes small. This coincides with the behavior of the graph of the mean boundary length in Figure 2(d) around the origin and the asymptotic formula

$$\mathbb{E}\mathcal{H}^1(\partial \Xi_r) \sim \frac{\pi T}{r \log^2 r}, \qquad r \to 0,$$

proved in [12].

The mean Euler number of $\Xi_r$, that is, 1 minus the mean number of holes, tends to $-\infty$ and its slope tends to infinity as $r \to 0$ since $\Xi_r$ is connected and the number of holes grows unboundedly. This behavior is reflected in Figure 2(f). For growing radius, the mean Euler number tends to 1 as more and more holes get filled by the dilation, and equals 1 when all holes disappear; see Figure 2(e).

Berezhkovskiĭ, Makhnovskiĭ and Suris [1] give an explicit formula for the expected area of the Wiener sausage in dimensions greater than 1. Its special case for $d = 2$ can be written as

$$\mathbb{E}V_{\Xi}(r) = \pi r^2 + \frac{8r^2}{\pi} \int_0^\infty \frac{1 - e^{-y^2 T/(2r^2)}}{y^3(J_0^2(y) + Y_0^2(y))} \, dy,$$

where $J_0$ and $Y_0$ are Bessel functions of order zero. Thus, $\mathbb{E}V_{\Xi}(r)$ and its derivative with respect to $r$—the expected boundary length $\mathbb{E}\mathcal{H}^1(\partial \Xi_r)$—can be computed numerically, in addition to the simulations shown above. Our computations show that the difference between the simulated and the approximated values is less than 10% for radii greater than 20 (resp., 200) for the boundary length (resp., area). For smaller radii, the difference gets bigger (up to 39%) due to the high slope of the boundary length and area, respectively.

7.2. *Approximation formulae.* The curves of the estimated mean intrinsic volumes were approximated with the *Curve Fitting Toolbox* in MATLAB using a trust region algorithm. The mean estimated area [see Figures 2(a) and 2(b)] of the Wiener sausage can be approximated by the function

$$\mathbb{E}V_{\Xi}(r) \approx \frac{441.83}{|\log r| + 265.265} + 387768 r^{0.252344} + 1064.51 r^{1.24465}$$
$$+ 0.235594 r^{2.27655}$$

depending on the dilation radius $r$. The absolute value of the relative error between this function and the estimated values is always less than 2% and the mean relative error is less than 0.6%.

For the boundary length in Figure 2(c), we fitted the function

$$\mathbb{E}\mathcal{H}^1(\partial \Xi_r) \approx \frac{54943.9}{r(\log r)^2 + 2.66290} + 93484.2 r^{-0.862668}$$
$$+ 9038.09 r^{-0.123680} + 15.6493 r^{0.901971}.$$

This is not exactly the derivative of the approximated mean area, but, nevertheless, the fit is very good. The relative errors are less than 1.9% for all radii and the mean relative error is less than 0.6%. In both cases, that is, for the estimated mean area and boundary length, their asymptotic behavior as $r \to 0$ coincides (up to a constant factor) with the formula of Le Gall given in the previous section.

A different class of functions was used to fit the Euler number. Here, a logarithmic term seems to be a good choice to approximate $2 - \mathbb{E}V_0(\Xi_r)$. [We did not approximate $\mathbb{E}V_0(\Xi_r)$ directly because it takes values close to zero in a certain interval and therefore one cannot compute relative errors to evaluate the fit.] Since this logarithmic function converges to 1 slowly as $r \to \infty$, we additionally used the cumulative distribution function $\Phi$ of the standard normal distribution to make it converge faster so that it is practically 1 for large radii. The resulting function fitted with MATLAB is given by

$$\mathbb{E}V_0(\Xi_r) \approx 1 - \frac{0.0423017(1 - \Phi((r - 224.899)/50.2096))}{\log(3.88182 \cdot 10^{-6} r^{1.88978} + 1.0) r^{0.153452}}.$$

The relative error of this fit is less than 10.5% and the mean relative error is less than 1.5%. For growing radii, the fitted function converges to one at an appropriate speed, which leads to low approximation errors. For radii greater than 160, the relative error is less than 4% and for radii greater than 560, the relative error is less than 1%.

**Acknowledgment.** The authors are grateful to Volker Schmidt for fruitful discussions.

J. RATAJ
CHARLES UNIVERSITY
FACULTY OF MATHEMATICS AND PHYSICS
SOKOLOVSKÁ 83
18675 PRAHA 8
CZECH REPUBLIC
E-MAIL: rataj@karlin.mff.cuni.cz

E. SPODAREV
D. MESCHENMOSER
ULM UNIVERSITY
INSTITUTE OF STOCHASTICS
HELMHOLTZSTR. 18
89069 ULM
GERMANY
E-MAIL: evgeny.spodarev@uni-ulm.de
       daniel.meschenmoser@uni-ulm.de